\documentclass[a4paper, reqno, 12pt]{amsart}

\usepackage{geometry}             

\usepackage{float}
\usepackage{bm}
\usepackage{fullpage,xcolor}
\usepackage[mathscr]{euscript}
\usepackage[all]{xy}
\usepackage{epsfig}
\usepackage{amsfonts}
\usepackage{mathptmx}

\usepackage{graphicx}
\usepackage{amssymb} 
\usepackage{amsmath}
\usepackage{amsthm}
\usepackage{mathrsfs}
\usepackage{epstopdf}
\usepackage{url}
\usepackage[msc-links,alphabetic]{amsrefs}
\usepackage{tikz}

\usepackage{color}
\makeatletter

\@addtoreset{equation}{section}
\makeatother 

\theoremstyle{plain}
\newtheorem{theorem}{Theorem}[section]
\newtheorem{corollary}[theorem]{Corollary}
\newtheorem{proposition}[theorem]{Proposition}

\newtheorem*{theorem3.14}{Theorem 3.14}
\newtheorem*{theorem3.15}{Theorem 3.15}

\theoremstyle{definition}
\newtheorem{definition}[theorem]{Definition}
\newtheorem{example}[theorem]{Example}
\newtheorem{remark}[theorem]{Remark}
\newtheorem{note}[theorem]{Note}

\usepackage{latexsym}
 
\title[Construction of weaving and polycatenane motifs \\ from periodic tilings of the plane]{Construction of weaving and polycatenane motifs \\ from periodic tilings of the plane}

\author{Mizuki Fukuda}
\address{Mathematics for Advanced Materials Open Innovation Laboratory, AIST, Tohoku University, 2-1-1 Katahira, Aoba-ku, Sendai 980-8577, Japan}
\email{mizuki.fukuda.d2@tohoku.ac.jp}

\author{Motoko Kotani}
\address{Advanced Institute for Materials Research, Tohoku University, 2-1-1 Katahira, Aoba-ku, Sendai 980-8577, Japan}
\email{motoko.kotani.d3@tohoku.ac.jp}

\author{Sonia Mahmoudi}
\address{Advanced Institute for Materials Research, Tohoku University, 2-1-1 Katahira, Aoba-ku, Sendai 980-8577, Japan;  RIKEN iTHEMS, 2-1 Hirosawa, Wako, Saitama 351-0198, Japan}
\email{sonia.mahmoudi@tohoku.ac.jp}

\subjclass[2020]{57K10, 57K12, 57K35, 57M10, 57M15, 05A05} 

\keywords{doubly periodic structures, DP tangles, weaves, polycatenanes, diagrams, motifs, links in the thickened torus, periodic tilings.}

\thanks{We would like to thank K. Shimokawa (Ochanomizu University) and M. Chas (Stony Brook University) for their precious comments and advice during this study. This work is supported by a Research Fellowship from JST CREST Grant Number JPMJCR17J4 and Grant-in Aid for JSPS Fellows Number 22J13397.}

\begin{document}

\begin{abstract} 
Doubly periodic \textit{(DP)} \textit{weaves} and \textit{polycatenanes} are complex entangled structures embedded in the Euclidean thickened plane, invariant under translations in two independent directions. Their topological properties are fully encoded within a quotient space under a periodic lattice, which we refer to as a \textit{motif}. On the diagrammatic level, a motif is a specific type of link diagram on the torus, consisting of essential closed curves for DP weaves or null-homotopic curves for DP polycatenanes. In this paper, we introduce a combinatorial methodology to construct these motifs from planar DP tilings using the concept of \textit{polygonal link transformations}. We also present an approach to predict the type of motif that can be constructed from a given DP tiling and a chosen polygonal link method. This approach has potential applications in various disciplines, such as materials science and chemistry.
\end{abstract}

\maketitle


\section{Introduction}\label{sec:1} 

Doubly periodic weaves (\textit{DP weaves}) and polycatenanes (\textit{DP polycatenanes}) are complex three-dimensional entangled structures consisting of non-intersecting curves, invariant under translations in two independent directions. They form two distinct classes of \textit{DP tangles}, defined in \cite{DLM} as the lift of links embedded in the thickened torus $T^2 \times I$. They can be differentiated by their type of components, namely open intertwined \textit{threads} for DP weaves and closed linked \textit{rings} for DP polycatenanes. These topological objects have applications in various fields, such as in materials science, biology, and chemistry, influencing the design of novel composites, molecular structures, and molecular assemblies, as mentioned in \cite{Evans:eo5019}, \cite{Evans:eo5020}, \cite{Evans:eo5050}, \cite{Hyde1}, \cite{Yaghi}, \cite{Panagiotou_2010}, \cite{ShimokawaBook}, \cite{Hyde2}.

In our previous work \cite{Sonia1}, we introduced a novel topological definition for specific classes of DP weaves, namely \textit{untwisted} and \textit{twisted weaves}. More generally, a DP weave can be described as an embedding of infinite simple open curves into the Euclidean thickened plane, denoted by $\mathbb{E}^2 \times I$, whose planar projection is a quadrivalent planar graph with crossing information at each vertex. In this paper, we present an equivalent definition for these DP weaves, using a knot theoretical approach. More specifically, we describe a DP weave as the lift to $\mathbb{E}^2 \times I$ of a link in $T^2 \times I$ consisting exclusively of essential closed curve components (see Definitions~\ref{def:untwistedweave} and~\ref{def:twistedweave}).

In the field of chemistry, polycatenanes are known for their structure consisting of multiple interlinked rings, with the simplest example being the Hopf link, which represents a classical pair of linked rings \cite{Yaghi}. More precisely, a polycatenane is an embedding of at least two linked rings. Polycatenanes have been also more generally defined as a disjoint union of knots in \cite{Shimokawa}. In the context of this paper, we consider each component of a DP polycatenane as a trivial knot, implying that each ring is an unknotted loop, as presented in Definition~\ref{def:polycatenane}.

As in classical knot theory (\cite{cromwell_2004}, \cite{DeToffoli2014-DETFAR}), we aim to study the topological properties of DP weaves and DP polycatenanes — and more generally of DP tangles (\cite{DLM}) — at the diagrammatic level. As initiated in \cite{Grishanov1}, \cite{Morton}, as well as \cite{KAWAUCHI2018230}, instead of analyzing a doubly periodic (DP) diagram on $\mathbb{E}^2$ containing an infinite number of crossings, we focus on its quotient under a point lattice isomorphic to $\mathbb{Z}^2$ acting on it by translations, namely a \textit{periodic lattice}. The resulting object is a link diagram on the torus $T^2$, which we refer to as a  \textit{motif}. A motif consists of essential (non-contractible) closed curves in the case of a \textit{weaving motif}, or null-homotopic (contractible) closed curves for a \textit{polycatenane motif}. A motif containing both essential and null-homotopic simple closed curves is referred to as a \textit{mixed motif}.

The primary aim of this paper is to construct weaving, polycatenane, and mixed motifs. Our strategy is to apply and formalize methods used to construct \textit{polyhedral links}, defined in \cite{Qiu1}, \cite{Qiu2}, and \cite{Qiu3}, to the context of periodic planar tilings \cite{TilingBook}. These methods involve transforming a polyhedron into a link by substituting each edge with a \textit{single} or \textit{double line}, possibly including \textit{twists}, and replacing each vertex with a set of \textit{branched} or \textit{crossed curves}. When applied to a planar tiling, we refer to these methods as \textit{polygonal link methods}. Note that the transformation of a periodic tiling into a weaving diagram by assigning over or under crossing information at each vertex has been explored in a separate context in \cite{ADAMS2019262} and \cite{ADAMS2020107045}.

\begin{figure}[ht]
\begin{center}
\includegraphics[width=5.3in]{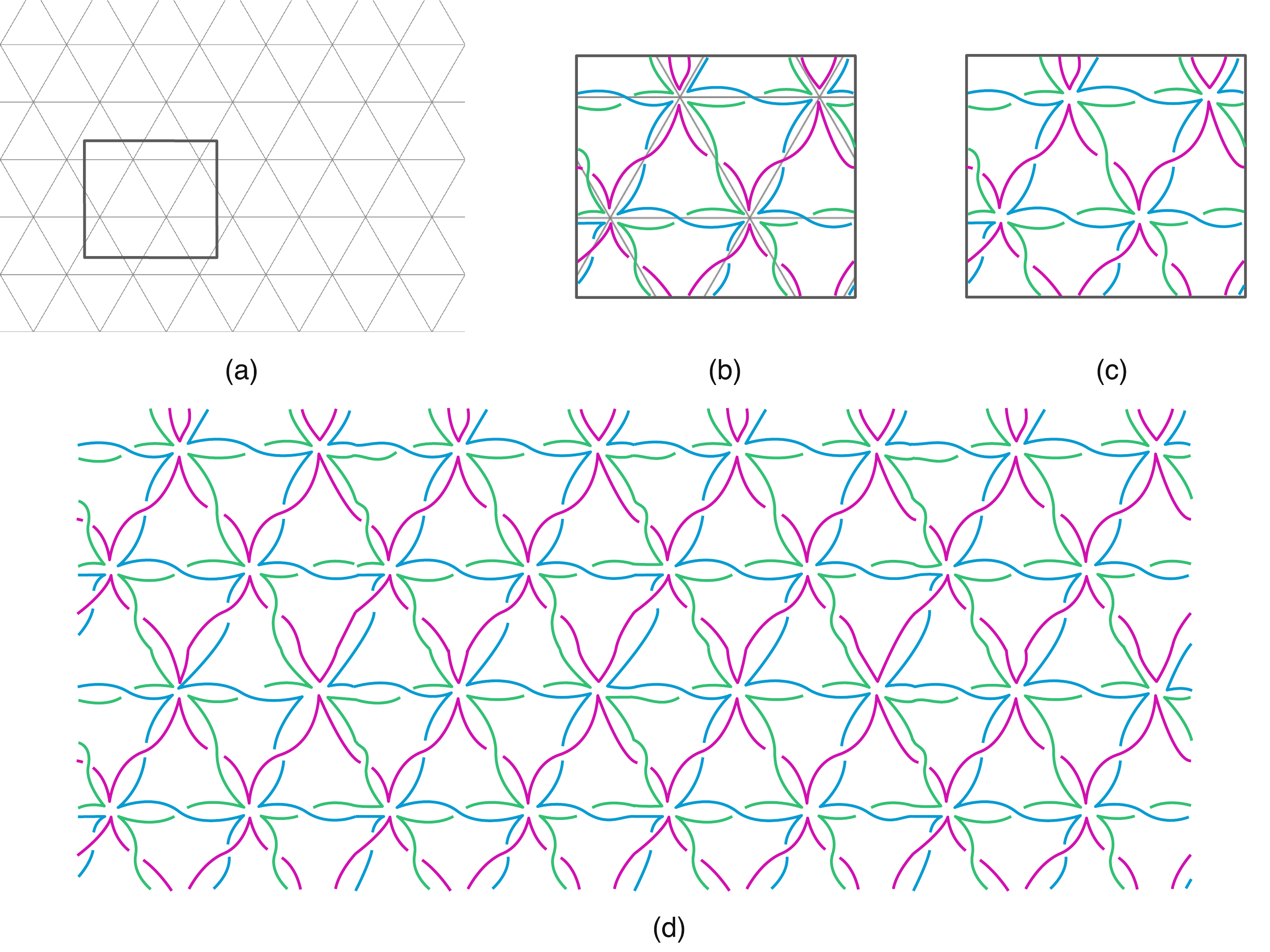}
\end{center}
\caption{\label{tiling-tcell} A DP tiling with a T-cell (a), transformed by a polygonal link method (b), resulting in a weaving motif (c), which lifts to a DP weave (d).}
\end{figure}

We construct motifs using the polygonal link methods as follows. First, let $\mathcal{T}$ be a doubly periodic (DP) tiling of the plane $\mathbb{E}^2$, and $\mathcal{P}$ be the quotient of $\mathcal{T}$ under a periodic lattice. This quotient is referred to as a $\mathcal{T}$\textit{-cell} and forms a connected graph on the (flat) torus $T^2$. Depending on the chosen polygonal link method, the vertices and edges of $\mathcal{P}$ are replaced with arcs and crossings, producing a motif as detailed in Section~\ref{sec:3-1}. An example is illustrated in Figure~\ref{tiling-tcell}. For simplicity, the polygonal link methods described in this paper construct \textit{alternating} motifs by convention. By an alternating motif, we mean that crossings alternate between over and under as one travels along each component of the motif \cite{MurasugiBook}. However, the prediction methods described in this paper do not rely on specific crossing information, so the motifs can be made non-alternating by adjusting the over and under crossing information.

This extension of the strategy of Qiu et al. to periodic planar graphs offers a new way to construct motifs of DP tangles. Using a graph-theoretic approach to study the properties of knots and links has proven useful in many cases, such as in \cite{PurcellAGT}, \cite{PurcellLondon}, \cite{ENDO20101002}, \cite{HydeEvans}, \cite{KOBATA2010213}. In the context of this paper, we are interested in predicting whether a chosen polygonal link method, denoted by $(\Lambda,L)$, applied to a $\mathcal{T}$-cell will create a weaving motif, a polycatenane motif, or a mixed motif. Each curve component $\gamma$ generated by $(\Lambda,L)$ forms a path in $T^2$, called a \textit{characteristic loop} and denoted  $\Delta^{\gamma}_{(\Lambda,L)}$, as defined in Section~\ref{sec:3-2}. By combining algebraic and combinatorial arguments, the second objective of this paper is to characterize the characteristic loops of a given $\mathcal{T}$-cell in order to predict the type of structure constructed in $\mathbb{E}^2 \times I$. Our main Theorems~\ref{thm:weavingprediction} and~\ref{thm:polycatenaneprediction} stated in Section~\ref{sec:3-3}, can find applications in many disciplines, particularly when considering a computational algorithm implementation of our results.

This paper is organized as follows: in Section~\ref{sec:2}, we present the definitions of weaving, polycatenane, and mixed motifs, as well as their corresponding DP tangles and diagrams. In Section~\ref{sec:3}, we formalize the concept of polygonal link methods and then prove our main prediction theorems.


\section{Preliminaries on DP weaves, DP polycatenanes and DP mix}\label{sec:2}

In this section, we present three classes of doubly periodic entangled structures (DP tangles) embedded in the Euclidean thickened plane, denoted as $\mathbb{E}^2 \times I$, where $\mathbb{E}^2$ represents the Euclidean plane and $I = [-1,1]$. 
To do so, let first $B = <u, v>$ be a basis of $\mathbb{E}^2$. We then consider the covering map $\rho: \mathbb{E}^2 \rightarrow T^2$ which assigns the longitude of the torus $T^2$ to $u$ and the meridian to $v$. Note that this covering map is extended to the product space by $\tilde{\rho}: \mathbb{E}^2 \times I \rightarrow T^2 \times I$.

\subsection{DP Weaves and their motifs}\label{sec:2-1}

\begin{figure}[ht]
\begin{center}
\includegraphics[width=5.5in]{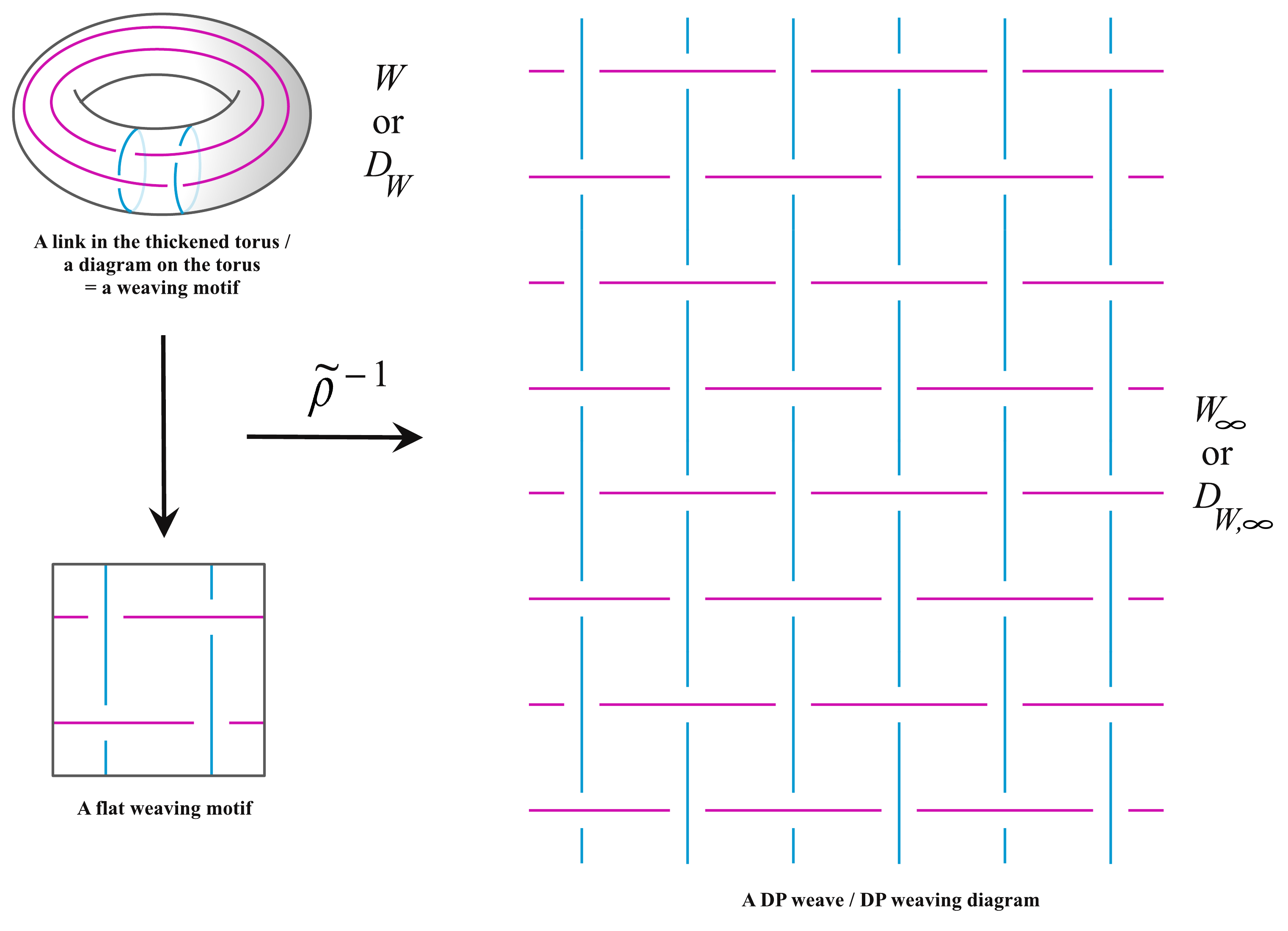}
\end{center}
\caption{\label{untwisted} 
An untwisted weaving motif (on the left), and its corresponding untwisted DP weave (diagram) (on the right).}
\end{figure}

In this subsection, we present the definitions of doubly periodic weaves (DP weaves) and their diagrams, as introduced in our previous work \cite{Sonia1}, to the context of this study. More specifically, we define a DP weave as the lift to $\mathbb{E}^2 \times I$ of a particular type of link embedded in $T^2 \times I$, characterized by its components being distinct \textit{torus knots}. Recall that a $(p,q)$-torus knot is an essential simple closed curve that intersects the torus meridian $p$ times and the torus longitude $q$ times, as detailed in \cite{MurasugiBook}(Chapter 7). Furthermore, a $(p',q')$-torus link consists of $g \geq 1$ isotopic components, each being a $(p,q)$-torus knot, where $p' = gp$ and $q' = gq$. In particular, the diagram of a torus link contains no self-intersections, and thus, we say that two components belonging to the same torus link (resp. their lifts under $\rho$) are \textit{parallel} in $T^2$ (resp. in $\mathbb{E}^2$). Moreover, we say that a single torus knot is parallel to itself in $T^2$ while its lifted curves in $\mathbb{E}^2$ are parallel to each other. Recall that a \textit{diagram} is defined as the canonical projection of a link in $T^2 \times I$ onto the torus $T^2$ via the map $\pi: T^2 \times I \to T^2$, where each double point includes over and under information, known as a \textit{crossing}. We can then define an untwisted DP weave as follows (see Figure~\ref{untwisted}).

\begin{definition}\label{def:untwistedweave} 
Let $W$ be a link embedded in $T^2 \times I$ with its corresponding diagram $D_W$ in $T^2$ satisfying the following conditions, with $i,j$ being positive integers, 
\begin{enumerate}
\item each component $t_i$ of $W$ is a $(p_i,q_i)$-torus knot, called \textit{thread}, where $p_i$ and $q_i$ are coprime integers,
 \smallbreak
 \item two parallel threads do not cross in $D_W$,
 \smallbreak
\item $W$ contains at least two distinct non-parallel threads $t_i$ and $t_j$.
\end{enumerate}
\noindent Then, the lift of $W$ (resp. $D_W$) under $\tilde{\rho}$ (resp. $\rho$) to $\mathbb{E}^2 \times I$ (resp. $\mathbb{E}^2 $) is called an \textit{untwisted DP weave}, denoted by $W_{\infty}$ (resp. an \textit{untwisted DP weaving diagram}, denoted by $D_{W, \infty}$). Moreover, $D_W$ is said to be an \textit{untwisted weaving motif} of $W_{\infty}$ (resp. $D_{W, \infty}$).
\end{definition}

\begin{figure}[ht]
\begin{center}
\includegraphics[width=5in]{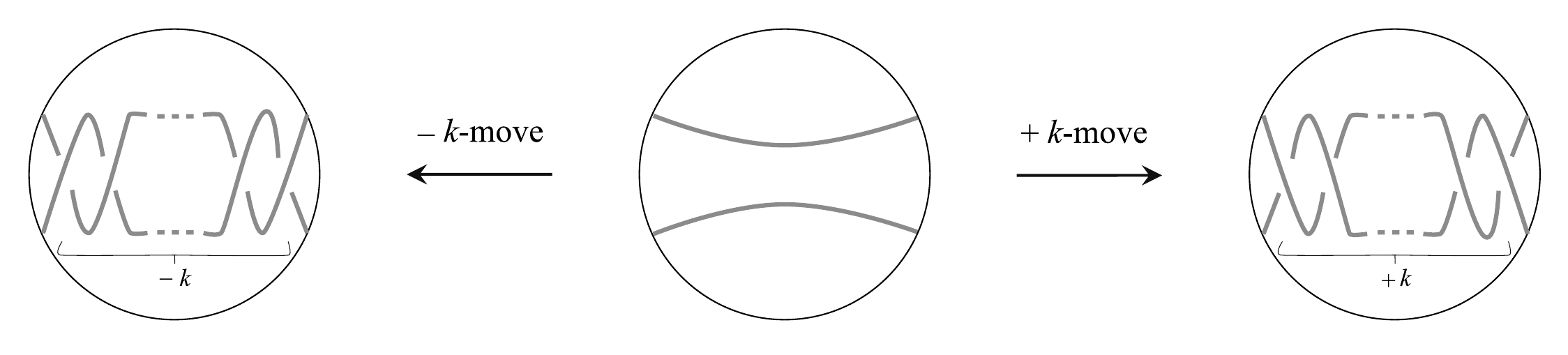}
\end{center}
\caption{\label{k-move} A $\pm k$-move (left and right) applied to a trivial tangle diagram (center).}
\end{figure}

Now, consider the following local surgery applied to an untwisted weaving motif $D_W$. Let $d$ be a topological disk such that the intersection $T = D_W \cap d$ forms a trivial tangle diagram 
containing exactly two arcs without crossings from the same thread $t_i$, or possibly from two parallel threads $t_i$ and $t_j$, for some $i$ and $j$ assuming the existence of $t_j$. By replacing the two parallel arcs in $T$ with $k$ half twists (see Figure~\ref{k-move} for an illustration), we introduce $k$ crossings into $D_W$. The resulting tangle diagram within $d$ is referred to as a \textit{twisted region}. This local operation is known as a $+k$-move for $k$ right-handed twists, or a $-k$-move for $k$ left-handed twists, as described in Chapter 3 of \cite{AdamsBook}.

\begin{figure}[ht]
\includegraphics[width=5.5in]{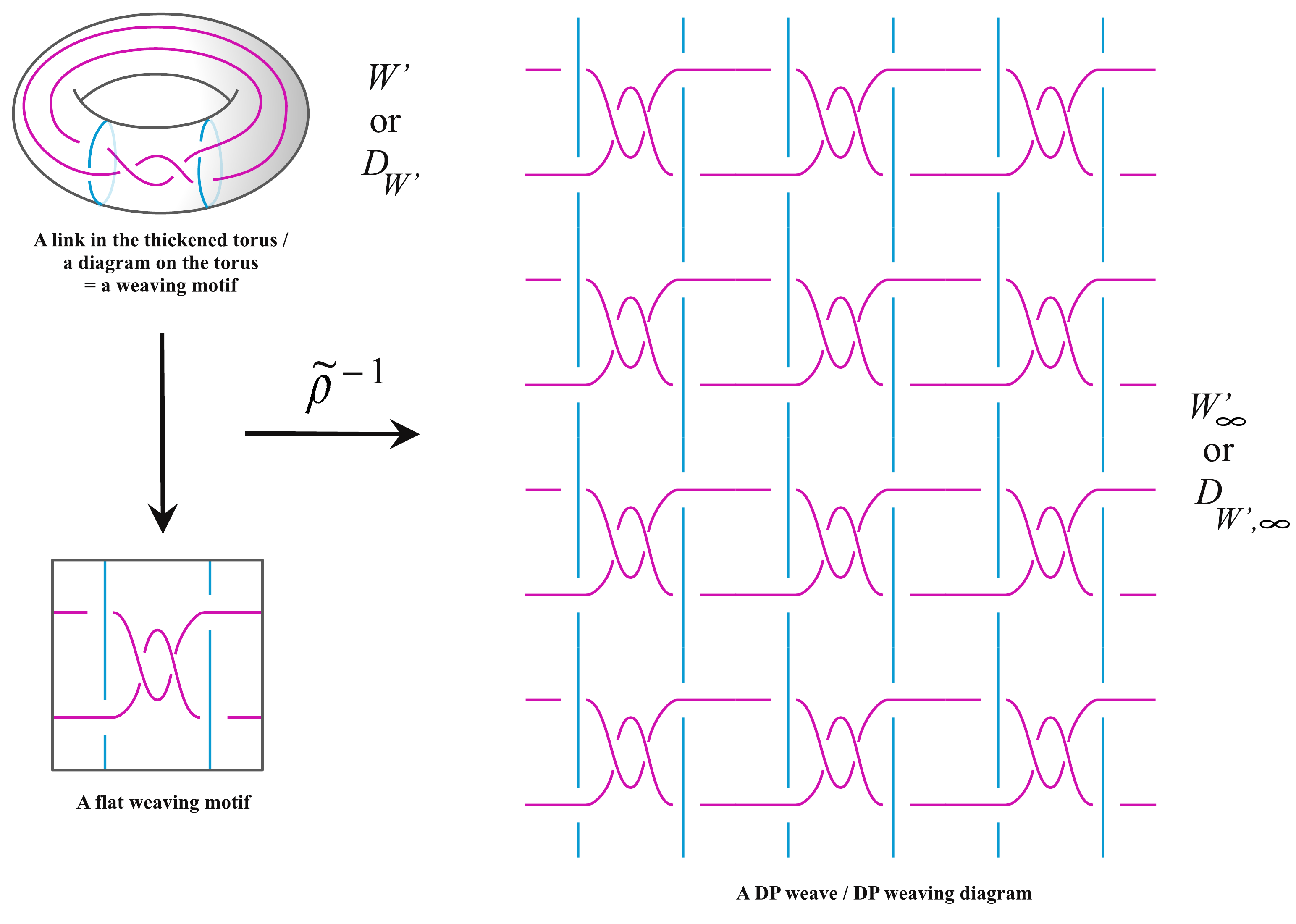}
\caption{\label{twisted} 
A twisted weaving motif (on the left), and its corresponding twisted DP weave (diagram) (on the right).}
\end{figure}

Note that parallel threads of $D_W$ (resp. $D_{W, \infty}$) are also said to be \textit{parallel} in $T^2$ (resp. in $\mathbb{E}^2$) after introducing twisted regions. The definition of a \textit{twisted DP weave} (see Figure~\ref{twisted}) follows from Definition~\ref{def:untwistedweave}.

\begin{definition}\label{def:twistedweave} 
Let $D_W$ be an untwisted weaving motif in $T^2$ where finitely many twisted regions are introduced. If the resulting diagram $D_{W'}$ on $T^2$ is not an untwisted weaving motif, then $D_{W'}$ is said to be a \textit{twisted weaving motif} and its lift $D_{W',\infty}$ to $\mathbb{E}^2$ under $\rho$ is referred to as a \textit{twisted DP weaving diagram}. Moreover, the lift of the corresponding link $W'$ in $T^2 \times I$ to $\mathbb{E}^2 \times I$ under $\tilde{\rho}$ is called a \textit{twisted DP weave} $W'_{\infty}$.
\end{definition}

\begin{figure}[ht]
\begin{center}
\includegraphics[width=5in]{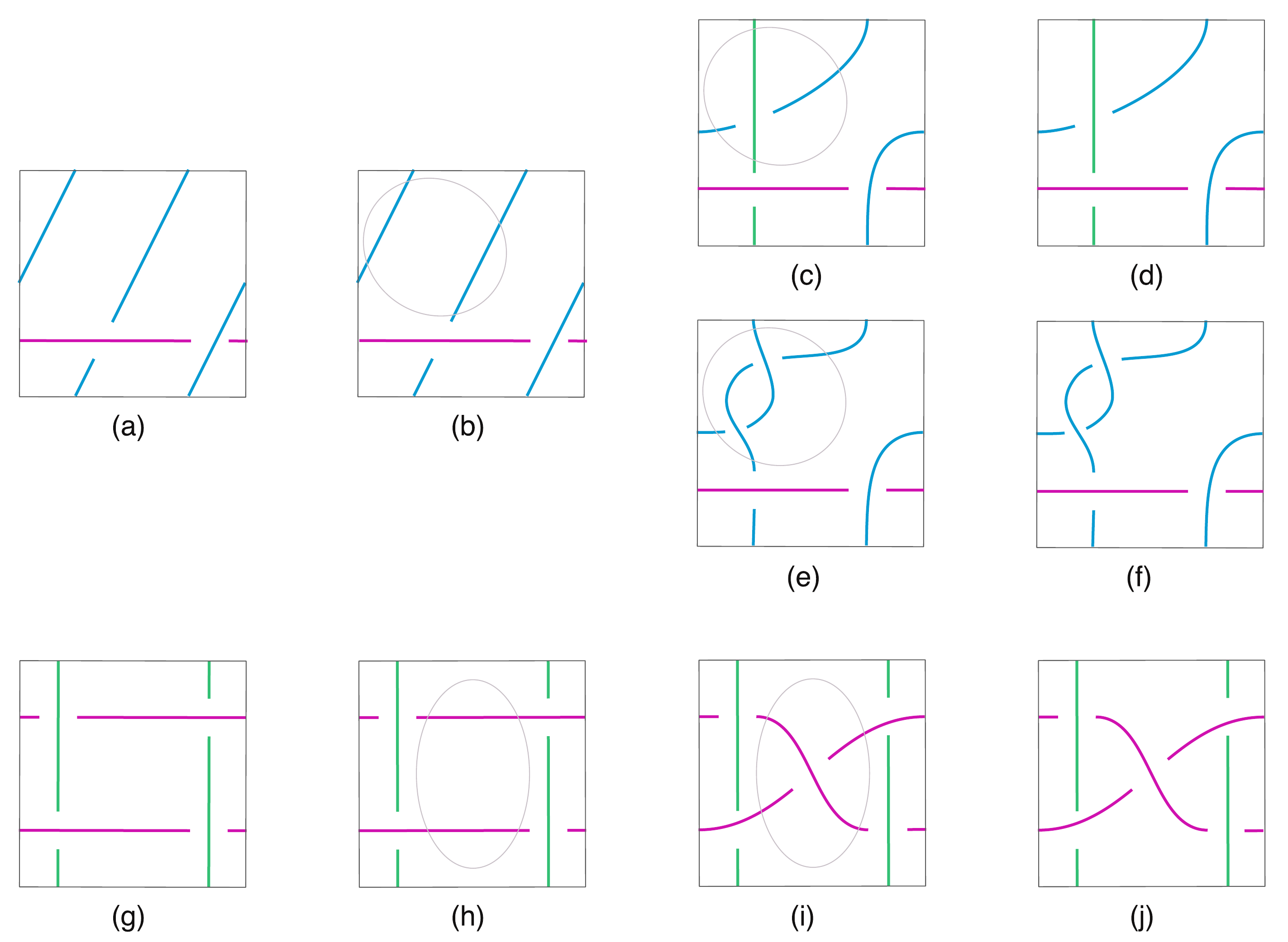}
\end{center}
\caption{\label{example-twisted} Examples on the effect of a $\pm k$-move on untwisted weaving motifs.}
\end{figure}

\begin{remark}\label{rem:weaving-motifs}
The introduction of finitely many twisted regions in an untwisted weaving motif can affect the resulting diagram in multiple ways. In particular, it can lead to an increase or decrease in the number of threads. A full investigation of these changes is outside the scope of this paper but may be of interest for future work. For example, consider an untwisted weaving motif $D_W^1$ (Figure~\ref{example-twisted} (a)) containing a $(2,1)$-torus knot component (in blue), and introduce a $\pm 1$-move (Figure~\ref{example-twisted} (c)) in a trivial tangle made of two arcs from this torus knot (Figure~\ref{example-twisted} (b)). The resulting weaving motif (Figure~\ref{example-twisted} (d)) -- whose number of components has increased by one (a blue $(1,1)$-torus knot, a pink $(1,0)$-torus knot and a green $(0,1)$-torus knot) -- remains untwisted, as does its lifted DP weaving diagram and the corresponding DP weave. 
However, if a $\pm 2$-move (Figure~\ref{example-twisted} (e)) is applied instead of the $\pm 1$-move, the resulting weaving motif (Figure~\ref{example-twisted} (f)) has the same number of components as $D_W^1$ but is now twisted. Now consider another example where an untwisted weaving motif $D_W^2$ (Figure~\ref{example-twisted} (g)) contains a $(2,0)$-torus link (in pink). By introducing a $\pm 1$-move (Figure~\ref{example-twisted} (i)) in a trivial tangle made of two arcs from this torus link (Figure~\ref{example-twisted} (h)) -- belonging to two different $(1,0)$-torus knots of the torus link -- the number of components in the diagram decreases by one, and the resulting motif becomes a twisted weaving motif (Figure~\ref{example-twisted} (j)).
\end{remark}

\smallbreak

\subsection{DP polycatenanes, DP mix and their motifs}\label{sec:2-2} 

\begin{figure}[ht]
\includegraphics[width=5.5in]{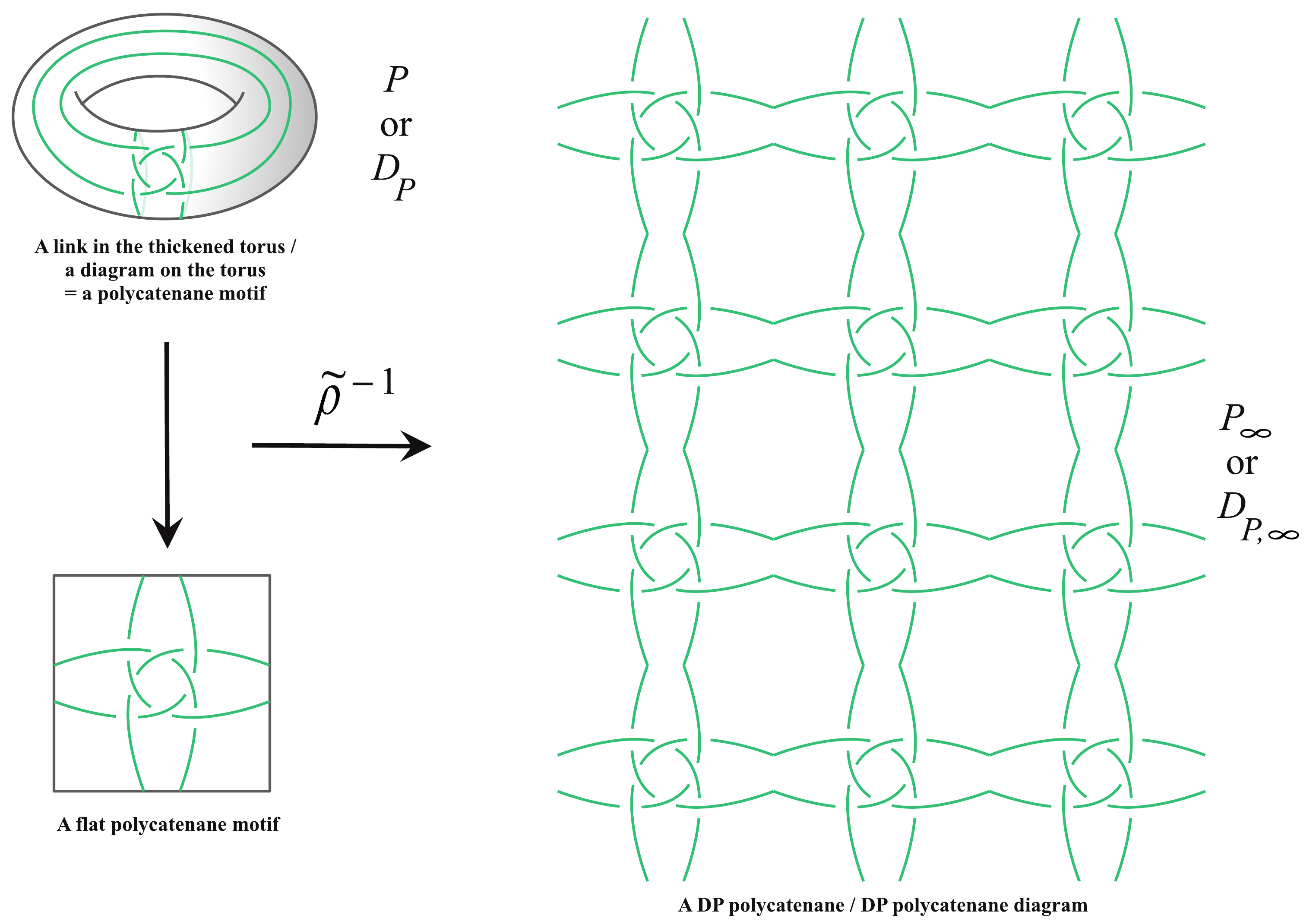}
\caption{\label{polycatenane} 
A polycatenane motif (on the left), and its corresponding DP polycatenane (diagram) (on the right).}
\end{figure}

In this subsection,  we extend the concept of DP weaves by introducing the definition of a DP polycatenane as the lift to $\mathbb{E}^2 \times I$ of a different type of link embedded in the thickened torus $T^2 \times I$ (see Figure~\ref{polycatenane} for an example). 

\begin{definition}\label{def:polycatenane} 
Let $P$ be a link embedded in $T^2 \times I$ with its corresponding diagram $D_P$ in $T^2$, such that each component of $P$ is a null-homotopic trivial knot, called \textit{ring}. Then the lift of $P$ (resp. $D_P$) under $\tilde{\rho}$ (resp. $\rho$) to $\mathbb{E}^2 \times I$ (resp. $\mathbb{E}^2 $) is called a \textit{DP polycatenane}, denoted by $P_{\infty}$ (resp. a \textit{DP polycatenane diagram}, denoted by $D_{P, \infty}$). Moreover, $D_P$ is said to be a \textit{polycatenane motif} of $P_{\infty}$ (resp. $D_{P, \infty}$).
\end{definition}

\begin{figure}[ht]
\includegraphics[width=5.5in]{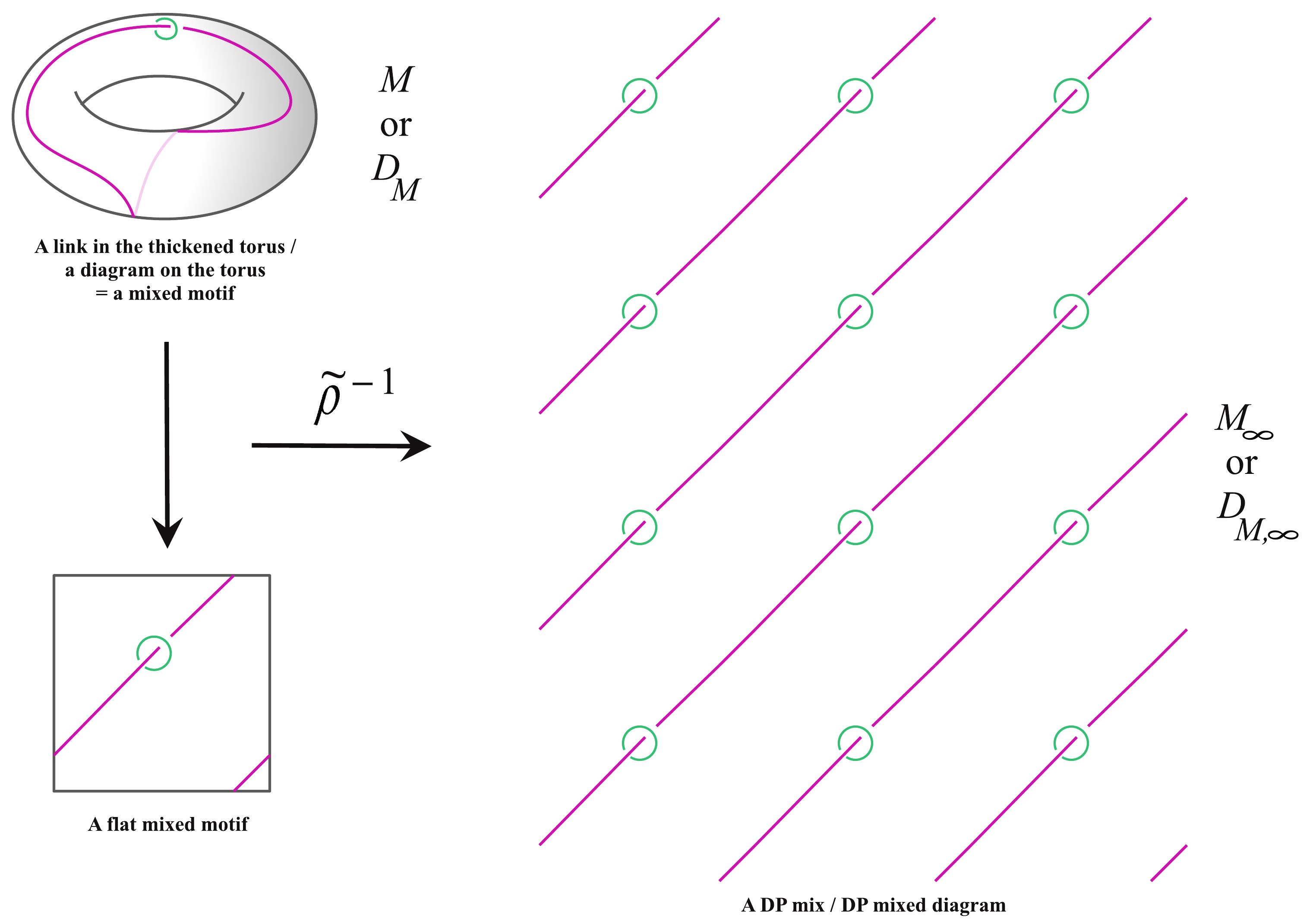}
\caption{\label{mixed} 
A mixed motif (on the left), and its corresponding  DP mix (diagram) (on the right).}
\end{figure}

Finally, we conclude this section by introducing a third class of DP tangles, that can be defined as a `mixed' structure containing both null-homotopic and essential components as follows (see Figure~\ref{mixed} for an illustration).

\begin{definition}\label{def:mixed} 
Let $M$ be a link embedded in $T^2 \times I$ with its corresponding diagram $D_M$ in $T^2$, such that the set of components of $M$ contains at least a thread and at least a ring. Then, the lift of $M$ (resp. $D_M$) under $\tilde{\rho}$ (resp. $\rho$) to $\mathbb{E}^2 \times I$ (resp. $\mathbb{E}^2 $) is called a \textit{DP mix}, denoted by $M_{\infty}$ (resp. a \textit{DP mixed diagram}, denoted by $D_{M, \infty}$). Moreover, $D_M$ is said to be a \textit{mixed motif} of $M_{\infty}$ (resp. $D_{M, \infty}$).
\end{definition}

\section{Motifs constructed by the polygonal link methods}\label{sec:3} 

In this section, we first provide a formal mathematical description of the polygonal link methods, each of which transforms a graph on the torus $T^2$ into a link diagram on $T^2$. Then, using algebraic and combinatorial arguments, we outline a strategy to predict the type of motif, namely a weaving, a polycatenane, or a mixed motif, that a given polygonal link method can create when applied to a selected $T$-cell.


\subsection{Polygonal link methods}\label{sec:3-1}

\begin{figure}[ht]
\includegraphics[width=5.5in]{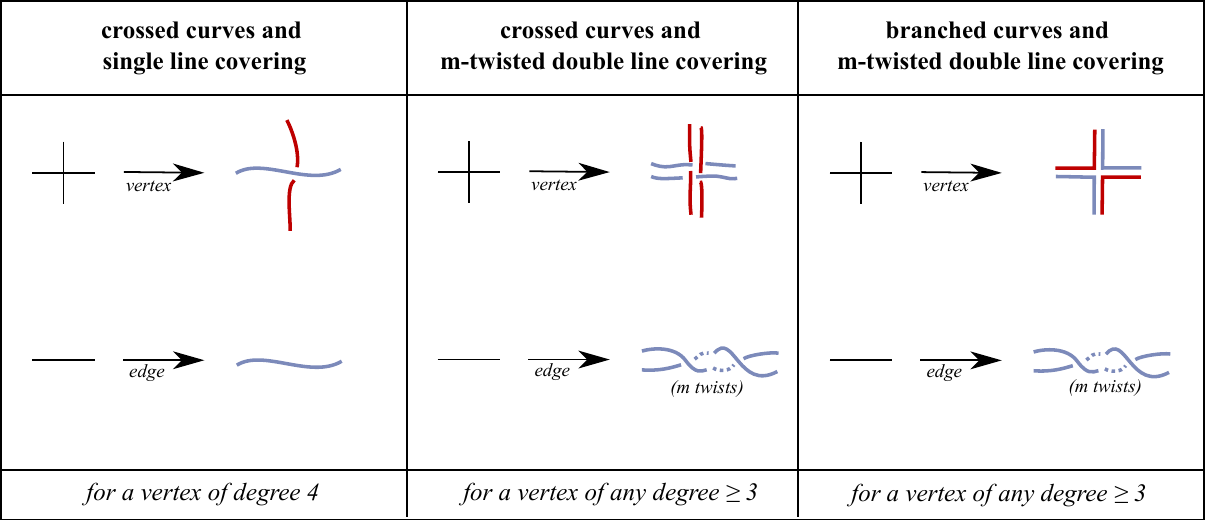}
\caption{\label{polygonal} The three polygonal link methods for a vertex and its adjacent edges.}
\end{figure}

In \cite{Qiu1}, \cite{Qiu2}, \cite{Qiu3}, new methods to transform a polyhedron into a `polyhedral link' are introduced, which consist of transforming each edge and each vertex of the polyhedron by one of the methods illustrated in Figure~\ref{polygonal}. In this subsection, we apply these transformations to graphs embedded on the (flat) torus and refer to them as \textit{polygonal link methods}.

We first recall that a \textit{polygonal tessellation} is a covering of the plane by polygons, called \textit{tiles}, such that the interior points of the tiles are pairwise disjoint. Moreover, two polygonal tessellations are said to be \textit{topologically equivalent}, if one can be mapped to the other by a homeomorphism from  $\mathbb{E}^2$ to itself, as defined in \cite{TilingBook}. We will refer to any such doubly periodic polygonal tessellation $\mathcal{T}$ of $\mathbb{E}^2$ as a \textit{DP tiling}.

Let now $\mathcal{P}$ be the quotient of the DP tiling $\mathcal{T}$ by a \textit{periodic lattice}, that is a group of discrete translational symmetry in two directions, isomorphic to $\mathbb{Z}^2$. We say that $\mathcal{P}$ is a \textit{$\mathcal{T}$-cell} of $\mathcal{T}$. Our strategy is to establish a methodical way to construct motifs by simultaneously transforming all the vertices and edges of $\mathcal{P}$ using the same polygonal link method. Here by an edge, we mean either a segment joining two distinct vertices or a loop that joins a vertex to itself. We introduce a formal definition of the three methods illustrated in Figure~\ref{polygonal}, considering a right-handed orientation of the plane.

\begin{definition}\label{def:3-1}
Let $\mathcal{P}$ be a $\mathcal{T}$-cell, and let $e_{v,v'}$ be the edge of $\mathcal{P}$ connecting the vertices $v$ and $v'$ in $\mathcal{P}$, possibly with $v=v'$.
Then $e_{v,v'}$ is said to be transformed by,
 \begin{enumerate} 
  \item a \textit{single line covering}, if $e_{v,v'}$ is replaced by a single strand $s_{v,v'}$, where both $v$ and $v'$ are vertices of degree four.
  \smallbreak
  \item an $\pm m$-\textit{twisted double line covering} if $e_{v,v'}$ is replaced by a pair of strands $s_{v,v'}$ and $s'_{v,v'}$ realizing a $\pm m$-move, where $m \geq 0$ is a positive integer.
\end{enumerate}
\end{definition}

\begin{figure}[ht]
\includegraphics[width=5in]{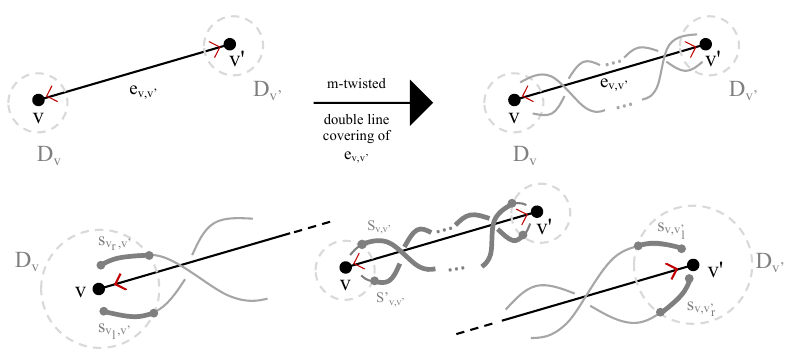}
\caption{\label{doubleline} Strand labeling for an $\pm m$-twisted double line covering.}
\end{figure}

\begin{remark}\label{rmk:leftright}
Note that if one of the vertices $v$ and $v'$ is not of degree four, then the single line covering transformation is not defined \cite{Qiu2}. Moreover, in the case of a $\pm m$-twisted double line covering, we decompose each strand into three connected segments:
$$ s_{v,v'} = s_{v_r,v'} \cup S_{v,v'} \cup s_{v,v'_\mu}, $$
$$ s'_{v,v'} = s'_{v_l,v'} \cup S'_{v,v'} \cup s'_{v,v'_\mu}, $$
with $\mu \in \{l,r\}$, where $l$ stands for \textit{left} and $r$ for \textit{right}. To assign the value $\mu$, consider a neighboring disk $d_v$ (with a right-handed orientation) of the vertex $v$ and orient the edge $e_{v,v'}$ toward $v$ in $d_v$. This implies that in a neighboring disk $d_{v'}$ (also with a right-handed orientation) of the vertex $v'$, the orientation of $e_{v,v'}$ is opposite since it is now oriented toward $v'$. Each such edge segment separates its neighborhoods into left and right regions, with respect to its orientation. Therefore, if an edge is covered by a $\pm m$-twisted double line, then in the interior of those disks, one of the strands is said to be on the left and the other on the right, as illustrated in Figure~\ref{doubleline}. We denote by $S_{v,v'}$ (resp. $S'_{v,v'}$) the curves that connect the strand on the right (resp. on the left) in $d_v$ to one of the strands in $d_{v'}$. 
We will use the following convention:
\smallbreak

 \begin{itemize} 
  \item If $m=0$ or $m$ is even, 
  $$ s_{v,v'} = s_{v_r,v'} \cup S_{v,v'} \cup s_{v,v'_l}, $$
  $$ s'_{v,v'} = s_{v_l,v'} \cup S'_{v,v'} \cup s_{v,v'_r}. $$
  \item If $m$ is odd, 
  $$ s_{v,v'} = s_{v_r,v'} \cup S_{v,v'} \cup s_{v,v'_r}, $$
  $$ s'_{v,v'} = s_{v_l,v'} \cup S'_{v,v'} \cup s_{v,v'_l}. $$
\end{itemize}
\end{remark}

\begin{figure}[ht]
\centering
   \includegraphics[width=3.5in]{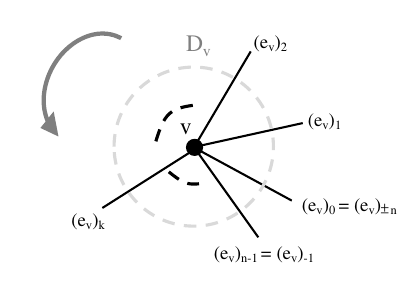}
      \caption{\label{edgelabel} Edge labeling at a neighborhood of a vertex $v$ of degree $n$.}
\end{figure}

Next, each vertex $v$ of degree $n \geq 3$ belonging to the given $\mathcal{T}$-cell will be transformed into a set of arcs with crossings, namely into a tangle diagram, in one of the three ways illustrated in Figure~\ref{polygonal} \cite{Qiu2}. In the neighborhood $d_v$ of $v$, we first label its adjacent edge segments starting from an arbitrary edge $e_v = (e_v)_0$. The first counterclockwise adjacent edge is $(e_v)_1$, the second is $(e_v)_2$, and so forth, until $(e_v)_n$, where $(e_v)_0 = (e_v)_n$. Conversely, by reading clockwise, the first adjacent edge is $(e_v)_{-1}$, the second is $(e_v)_{-2}$, and so on, up to $(e_v)_{-n}$, where $(e_v)_0 = (e_v)_{-n}$. See Figure~\ref{edgelabel} for an illustration.

\bigbreak

Now, consider a single line covering transformation applied to an edge $e_{v,v'}$. This transformation implies that both adjacent vertices $v$ and $v'$ of $e_{v,v'}$ are of degree four and are replaced by a \textit{4-crossed curves} transformation, since the single line covering transformation is not defined for vertices of other degrees in the polygonal link methods \cite{Qiu2}. In this case, each pair of opposite strands is glued to form a single strand. 

We use the convention that during the construction, the strand being the union of $(s_v)_1$ and $(s_v)_3$ will be positioned over the second strand, where a strand $(s_v)_i$ replaces the edge segment $(e_v)_i$ in $d_v$, for $i = 0, 1, 2, 3$. Such over-crossing information can be arbitrarily reversed into under-crossing information.

However, if an edge undergoes a $\pm m$-twisted double line covering transformation, two options arise for transforming its adjacent vertices according to the polygonal link methods \cite{Qiu1, Qiu2}, as exemplified in Figure~\ref{polygonal}. In one case, we can connect each right strand $(s_{v_r,v'})_i$ in $d_v$ with its closest neighboring left strand $(s_{v_l,v''})_{i+1}$, constructed from the edge adjacent to vertices $v$ and $v'' \neq v'$, by performing a \textit{branched curves} transformation \cite{Qiu1}. 

Alternatively, a \textit{crossed curves} transformation can be applied by connecting each right strand $(s_{v_r,v'})_i$ to the second counterclockwise adjacent left strand $(s_{v_l,v''})_{i+2}$ \cite{Qiu2}. For both transformations, we adopt the convention that in the resulting crossings, the left strand appears over the right strand. Yet, the transformation, independent of the crossing information, can be arbitrarily reversed.

\begin{definition}\label{def:vertex}
A vertex $v$ of degree $n_v \geq 3$ of $\mathcal{P}$ is said to be transformed by,
 \begin{enumerate} 
  \item [a.] \textit{4-crossed curves}, if $n_v = 4$, and if each strand $(s_{v,v'})_i$ is glued to the strand $(s_{v,v''})_{i+2}$ to form a single strand. This occurs when each edge adjacent to $v$ has been replaced by a single line, and $v, v', v''$ are distinct adjacent vertices.
  \smallbreak
  \item [b.] \textit{crossed curves}, if each strand $(s_{v_r,v'})_i$ connects with the strand $(s_{v_l,v''})_{i+2}$ to form a single strand. This occurs when each edge adjacent to $v$ has been replaced by an $\pm m$-twisted double line, where $m \geq 0$, and $v, v',v ''$ are distinct adjacent vertices.
  \smallbreak
  \item [c.] \textit{branched curves}, if each strand $(s_{v_r,v'})_i$ connects with the strand $(s_{v_l,v''})_{i+1}$ to form a single strand. This occurs when each edge adjacent to $v$ has been replaced by an $\pm m$-twisted double line, where $m \geq 0$, and $v, v', v''$ are distinct adjacent vertices.
\end{enumerate}
\end{definition}

\begin{figure}[ht]
\centering
   \includegraphics[width=6in]{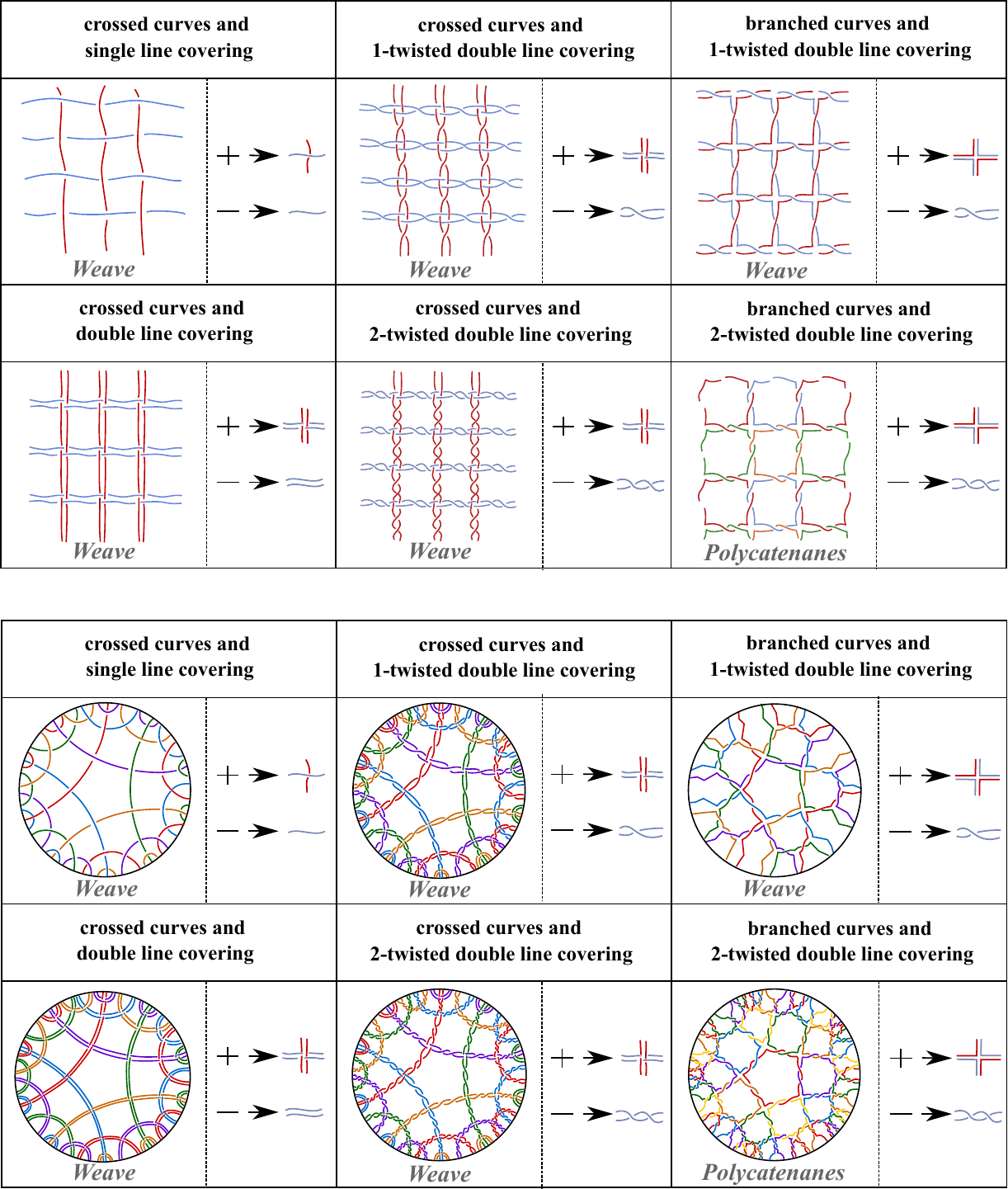}
      \caption{\label{examples} Examples of weaving and polycatenane diagrams in $\mathbb{E}^2$ (on the top) and $\mathbb{H}^2$ (on the bottom).}
\end{figure}

The transformations applied to the vertices and edges of a $\mathcal{T}$-cell can now be merged into three distinct polygonal link methods, as defined below. This definition formalizes an adaptation of the three polyhedral link methods presented in \cite{Qiu2} for graphs on the torus, corresponding to planar DP tilings. For illustrations, we refer to Figure~\ref{examples}.

\begin{definition}\label{def:polygonal}
A $\mathcal{T}$-cell is said to be transformed by the \textit{polygonal link method} $(\Lambda, L)$, where $L \in \{s, \pm m\}$ (with $m$ as a positive integer) and $\Lambda \in \{\mathit{C}_r, \mathit{B}_r\}$, if all its vertices and edges are transformed by the same method with,

$$
 (\Lambda,L) = \left\{
    \begin{array}{lll}
        (\mathit{C}_r,s) & \mbox{: 4-crossed curves and single line covering},  \\
        (\mathit{C}_r,\pm m) & \mbox{: crossed curves and $m$-twisted double line covering},  \\
        (\mathit{B}_r,\pm m) & \mbox{: branched curves and $m$-twisted double line covering}. \\
    \end{array}
\right.
$$
\end{definition}

\bigbreak

\begin{note}\label{rmk:combination}
   A generalization of these construction methods can be considered. For example, if $\mathcal{P}$ is a $\mathcal{T}$-cell containing at least two vertices, such as $v$ and $v'$, one might transform $v$ into a set of crossed curves and $v'$ into a set of branched curves, assuming a double line covering for each edge of $\mathcal{P}$. Similarly, for two distinct edges, $e$ and $e'$ of $\mathcal{P}$, $e$ can be replaced with a $\pm m$-twisted double line, and $e'$ with a $\pm m'$-twisted double line, where $m$ and $m'$ are possibly distinct. In subsequent subsections, we do not explore the scenario where two vertices undergo transformations using different methods. However, our main theorems remain valid if two edges are covered by different twisted regions, such as a $\pm k$-move and a $\pm k'$-move, provided that the parity of $k$ and $k'$ is the same.
\end{note}

\begin{remark}\label{rmk:3-5}
We believe that transformations of the edges and vertices of a $\mathcal{T}$-cell other than those summarized in Definition~\ref{def:polygonal} may lead to additional construction methods for motifs of DP tangles. However, these are beyond the scope of this paper and present interesting directions for further exploration. Moreover, note that the three polygonal link methods presented in this section can also be applied to periodic tilings of the hyperbolic plane $\mathbb{H}^2$. This application would lead to the construction of hyperbolic motifs on higher genus surfaces and to diagrams of periodic tangles in $\mathbb{H}^2$. In Figure~\ref{examples}, we use the Poincaré disk model for an illustration.
\end{remark}

\smallbreak
\subsection{Characteristic loops of polygonal link methods}\label{sec:3-2}

The second main objective of this paper is to predict whether a polygonal link method $(\Lambda, L)$ applied to a chosen $\mathcal{T}$-cell will generate a weaving motif, a polycatenane motif, or a mixed motif. More precisely, given a $\mathcal{T}$-cell with labeled and oriented edges, our approach is to characterize combinatorially the types of \textit{polygonal chains} that will be formed by the curve components created from $(\Lambda, L)$.

\begin{definition}\label{def:characteristicpath}
Let $\mathcal{P}$ be a $\mathcal{T}$-cell with oriented edges $\{e_1, \cdots, e_n\}$, and consider the polygonal link method $(\Lambda,L)$ applied to $\mathcal{P}$. 
Let also $\gamma$ be an oriented curve on $T^2$ defined by the finite union of connected strands covering the edges of $\mathcal{P}$ by $(\Lambda,L)$, starting arbitrarily from an edge $\delta_0$ and ending on the edge $\delta_j$, with $\delta_0, \, \delta_j \in \{e_1, \cdots, e_n\}$. Then, $\gamma$ is characterized by the ordered and reduced sequence $\Delta^{\gamma}_{(\Lambda,L)}$ of adjacent edges of $\mathcal{P}$, such that,
\begin{equation} 
\Delta^{\gamma}_{(\Lambda,L)} = (\delta_0^{\pm}, \delta_1^{\pm}, \cdots , \delta_j^{\pm})_{(\Lambda,L)}, 
\end{equation}
where, for all $i \in \{1, \cdots, n\}$, $\delta_i \in \{e_1, \cdots, e_n\}$, and $\delta_i^+$ and $\delta_i^-$ denote the same edge with opposite orientations.
Then, the polygonal chain $\Delta^{\gamma}_{(\Lambda,L)}$ is called a \textit{characteristic path of} $\mathcal{P}$ for $(\Lambda,L)$.
\end{definition}

\begin{remark}
    In Definition~\ref{def:characteristicpath}, the term \textit{reduced} means that $\Delta^{\gamma}_{(\Lambda,L)}$ does not contain a repetition of the entire sequence from $\delta_0$ to $\delta_j$, although partial sequences and loops are allowed. 
\end{remark}

\smallbreak

In particular, note that a characteristic path of $(\Lambda,L)$ forms a loop on $\mathcal{P}$ as follows.

\begin{proposition}
    Let $\mathcal{P}$ be a $\mathcal{T}$-cell and $\Delta^{\gamma}_{(\Lambda,L)}$ be a characteristic path of $\mathcal{P}$ for $(\Lambda,L)$.
    Then, $\Delta^{\gamma}_{(\Lambda,L)}$ is a closed path on $\mathcal{P}$ and is called a \textit{characteristic loop of} $\mathcal{P}$ for $(\Lambda,L)$.
\end{proposition}

\bigbreak

\begin{proof}
A $\mathcal{T}$-cell $\mathcal{P}$ can be defined as a graph with a finite number of vertices and edges on $T^2$. Moreover, since $T^2$ is a compact surface, the path $\Delta^{\gamma}_{(\Lambda,L)}$ on $\mathcal{P}$ cannot escape into infinity and is bounded by definition. By Definition~\ref{def:characteristicpath}, since each edge has a specific orientation and can only connect two specific vertices, any path $\Delta^{\gamma}_{(\Lambda,L)}$ must eventually use a previously visited vertex, as it can only stop if it starts repeating the full sequence of edges. Then, because the set of possible paths $\Delta^{\gamma}_{(\Lambda,L)}$ on $\mathcal{P}$ is finite, a path cannot infinitely extend without repeating. This cyclic property is guaranteed since any non-repeating extension of the path would require either an increase in the number of vertices or edges, or an escape into an unbounded region, none of which are possible on $T^2$.
\end{proof}

\begin{figure}[ht]
\centering
   \includegraphics[width=5.5in]{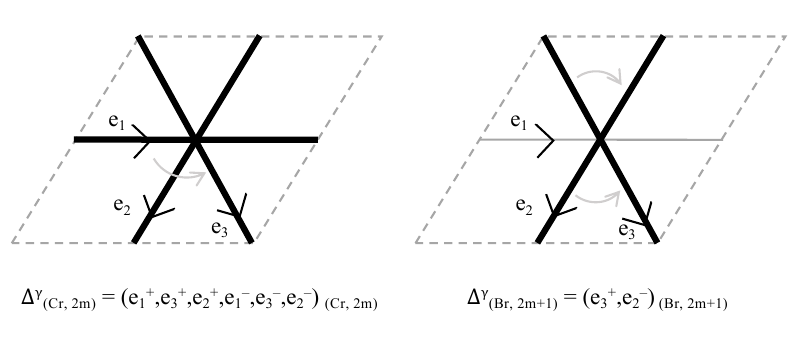}
      \caption{\label{triangleT} A $\mathcal{T}$-cell of a triangular tiling. On the left (resp. right), a characteristic loop corresponding to Figure~\ref{triangleU} (left) (resp. top). }
      \end{figure}

\begin{figure}[ht]
\centering
   \includegraphics[width=5in]{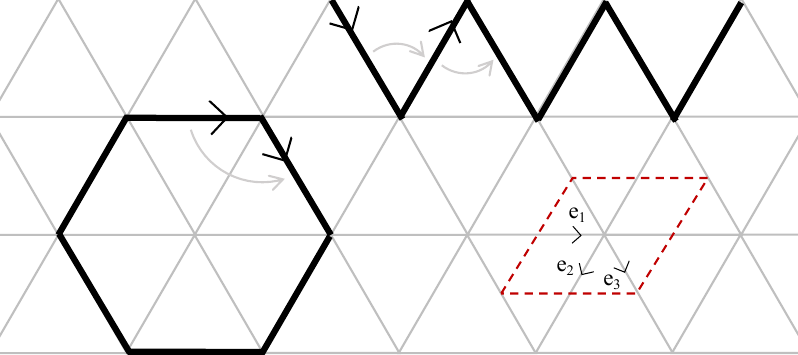}
      \caption{\label{triangleU} Triangular tiling.
      On the left, a simple closed polygonal chain that would be covered by applying $(\Lambda,L) = (C_r,\pm 2m)$.
      On the top, a subset of a simple open polygonal chain that would be covered by applying $(\Lambda,L) = (B_r,\pm 2m+1)$. 
      On the right, a $\mathcal{T}$-cell with labeled vertices and edges.}
\end{figure}

Figure~\ref{triangleT} and Figure~\ref{triangleU} illustrate two examples for a triangular tiling. In Figure~\ref{triangleT}, two characteristic loops are shown for two different polygonal link methods, and their lifts to the plane are (partially) presented in Figure~\ref{triangleU}.

We now distinguish the characteristic loops for the three polygonal link methods.

\subsubsection{{\bf Polygonal link method: 4-Crossed curves and single line covering.}} \label{sec:3.2.1}

We first present the case $(\Lambda,L) = (\mathit{C}_r,s)$, where each vertex of a $\mathcal{T}$-cell has degree four and is transformed by 4-crossed curves, while each edge is covered by a single line.

Let $\mathcal{P}$ be a $\mathcal{T}$-cell with oriented edges and select an arbitrary vertex of $\mathcal{P}$ from which we start a characteristic loop $\Delta^{\gamma}_{(\mathit{C}r,s)}$. Then, by Definition~\ref{def:vertex}, each element of $\Delta^{\gamma}_{(\mathit{C}_r,s)}$ is the second counterclockwise adjacent edge to the previous one at their common vertex. Thus, by using the notations of Definition~\ref{def:characteristicpath} and Figure~\ref{edgelabel}, we have, 

$$\Delta^{\gamma}_{(C_r,s)} = (\delta_0^{\pm}, \delta_1^{\pm}, \cdots , \delta_k^{\pm})_{(C_r,s)},$$ 
where for all $i \in \{1, \cdots, k\}$
\begin{equation}
\delta_i = (\delta_{i-1})_2 \mbox{ and } \delta_0 = (\delta_{k})_2,
\end{equation}
where the index $2$ in $\delta_i = (\delta_{i-1})_2$ means that $\delta_i$ is the second counterclockwise edge adjacent to $\delta_{i-1}$.

\smallbreak

\subsubsection{{\bf Crossed curves and $\pm m$-twisted double line covering.}}\label{sec:3.2.2}

Now, we present the case $(\Lambda,L) = (\mathit{C}_r,\pm m)$, where each vertex of a $\mathcal{T}$-cell is transformed by crossed curves, while each edge is covered by a double line.

With the same notations, the strategy is similar to the previous case. The main difference concerns the $\pm m$-twisted double line covering and depends on the parity of $m$. By Definition~\ref{def:vertex} and using the above notations, along with Remark~\ref{rmk:leftright}, we can state that,

\begin{itemize}  
  \item [-] if $m$ is even, then each curve $\gamma$ constructed by the method $(\mathit{C}_r,\pm m)$ in $\mathcal{P}$ is an alternating union of a right strand, followed by a left strand, and so forth. Thus, $\gamma$ is covering a characteristic loop made of edges that are consecutive second counterclockwise adjacent edges at each vertex they cross.
 \smallbreak
  \item [-] if $m$ is odd, then each curve $\gamma$ constructed by the method $(\mathit{C}_r,\pm m)$ in $\mathcal{P}$ is an alternating union of two consecutive right strands, followed by two consecutive left strands, and so forth. In this case, the consecutive edges forming the characteristic loop alternate between second counterclockwise adjacent edge and second clockwise adjacent edge from the previous one.
\end{itemize}
Thus, using Definition~\ref{def:characteristicpath} and the notations of Figure~\ref{edgelabel}, we have,
$$\Delta^{\gamma}_{(C_r,\pm m)} = (\delta_0^{\pm}, \delta_1^{\pm}, \cdots , \delta_k^{\pm})_{(C_r,\pm m)},$$ 
where for all $i \in \{1, \cdots, k\}$, 
if $m=0$ or if $m$ is an even integer,
\begin{equation}
\delta_i = (\delta_{i-1})_2 \mbox{ and } \delta_0 = (\delta_{k})_2,
\end{equation}
where the index $2$ in $\delta_i = (\delta_{i-1})_2$ means that $\delta_i$ is the second counterclockwise edge adjacent to $\delta_{i-1}$.
Otherwise,
\begin{equation}
\delta_{\,2i} = (\delta_{\,2i-1})_2  \ \mbox{ , } \ \delta_{\,2i+1} = (\delta_{\,2i})_{-2} \mbox{ and } \delta_0 = (\delta_{k})_2,
\end{equation}
where the index $2$ (resp. $-2$) in $\delta_{\,2i} = (\delta_{\,2i-1})_2$ (resp. $\delta_{\,2i+1} = (\delta_{\,2i})_{-2}$) means that $\delta_{\,2i}$ (resp. $\delta_{\,2i+1}$) is the second counterclockwise (resp. clockwise) edge adjacent to $\delta_{\,2i-1}$ (resp. $\delta_{\,2i}$).

\subsubsection{{\bf Branched curves and $\pm m$-twisted double line covering.}}\label{sec:3.2.3}

Finally, we present the case $(\Lambda,L) = (\mathit{B}_r,\pm m)$, where each vertex of a $\mathcal{T}$-cell is transformed by branched curves, while each edge is covered by a double line. This follows directly from the case $(\Lambda,L) = (\mathit{C}_r,\pm m)$ by replacing the second (counter-)clockwise adjacent edge by the first one: 
$$\Delta^{\gamma}_{(B_r,\pm m)} = (\delta_0^{\pm}, \delta_1^{\pm}, \cdots , \delta_k^{\pm})_{(B_r,\pm m)},$$ 
where for all $i \in \{1, \cdots, k\}$, 
if $m=0$ or if $m$ is an even integer,
\begin{equation}
\delta_i = (\delta_{i-1})_1 \mbox{ and } \delta_0 = (\delta_{k})_1,
\end{equation}
where the index $1$ in $\delta_i = (\delta_{i-1})_1$ means that $\delta_i$ is the first counterclockwise edge adjacent to $\delta_{i-1}$.
Otherwise,
\begin{equation}
\delta_{\,2i} = (\delta_{\,2i-1})_1  \ \mbox{ , } \ \delta_{\,2i+1} = (\delta_{\,2i})_{-1} \mbox{ and } \delta_0 = (\delta_{k})_1,
\end{equation}
where the index $1$ (resp. $-1$) in $\delta_{\,2i} = (\delta_{\,2i-1})_1$ (resp. $\delta_{\,2i+1} = (\delta_{\,2i})_{-1}$) means that $\delta_{\,2i}$ (resp. $\delta_{\,2i+1}$) is the first counterclockwise (resp. clockwise) edge adjacent to $\delta_{\,2i-1}$ (resp. $\delta_{\,2i}$).

\smallbreak
\subsection{Weaving, polycatenane or mixed motifs}\label{sec:3-3}

Let $\mathcal{P}$ be a $\mathcal{T}$-cell of a DP tiling $\mathcal{T}$ and let $(\Lambda, L)$ be a polygonal link method. Since each curve component of a motif covers a characteristic loop in $\mathcal{P}$, we can identify its type using the definitions of weaving (Definitions~\ref{def:untwistedweave} and~\ref{def:twistedweave}), polycatenane (Definition~\ref{def:polycatenane}), and mixed motifs (Definition~\ref{def:mixed}). 

To do so, we first introduce the notion of `divided curves' as follows. On the diagrammatic level, recall that any closed curve component $\omega$ in $T^2$ of a motif is either essential or null-homotopic, and lifts to simple curves in $\mathbb{E}^2$. Moreover, $\omega$ is said to be \textit{simple} if it does not admit any self-intersections, or \textit{non-simple} if it admits a finite number of self-intersections, each being a transverse double point. Since both simple and non-simple closed curves in $T^2$ may lift to simple curves in $\mathbb{E}^2$, these two cases are considered in the study of the three types of motifs. However, we aim to exclude the case where a non-simple closed curve in $T^2$ lifts to self-intersecting curves in $\mathbb{E}^2$. 

In particular, if $\omega$ is a non-simple closed curve in $T^2$ with a self-intersection on a point $p$, then $p$ can divide $\omega$ into two closed curves, which may also be self-intersecting. We refer to these two closed curves in $T^2$ as $\omega_p$ and $\omega'_p$, which we call the \textit{divided curves} of $\omega$ for $p$. Notably, if at least one of $\omega_p$ or $\omega'_p$ is null-homotopic in $T^2$, then this divided curve is said to be \textit{trivial} and $\omega$ lifts to self-intersecting curves in $\mathbb{E}^2$. For a detailed justification, see \cite{RotationNumber}. 

We are now in a position to characterize the characteristic loops for each type of motif constructed from a given $\mathcal{T}$-cell and a chosen polygonal link method as follows.

\begin{proposition}\label{prop:characteristic} 
A pair $\big(\mathcal{P},(\Lambda,L) \big)$, where $\mathcal{P}$ is a given $\mathcal{T}$-cell, and $(\Lambda,L)$ is a polygonal link method, generates: 
\begin{enumerate}  
  \item a weaving motif if and only if all the characteristic loops in $\mathcal{P}$ are homotopic to essential closed curves, without trivial divided curves, and such that at least two of them are not parallels.   
  \smallbreak
  \item a polycatenane motif if and only if all the characteristic loops in $\mathcal{P}$ are null-homotopic closed curves, without trivial divided curves.
  \smallbreak
  \item a mixed motif if and only if the set of all the characteristic loops in $\mathcal{P}$ contains at least a null-homotopic and an essential closed curves, without trivial divided curves.  
\end{enumerate}
\end{proposition}

\begin{proof} 
To prove the forward implications of this proposition, we refer to the definitions of untwisted and twisted weaving motifs, polycatenane motifs and mixed motifs, namely Definitions~\ref{def:untwistedweave},~\ref{def:twistedweave},~\ref{def:polycatenane},~\ref{def:mixed}, respectively. 
Let $\mathcal{P}$ be a given $\mathcal{T}$-cell, and $(\Lambda,L)$ be a polygonal link method. If $\big(\mathcal{P},(\Lambda,L) \big)$ creates a weaving motif, then since all its component are essential curves in $T^2$, each characteristic loop of $\mathcal{P}$ for $(\Lambda,L)$ must be homotopic to an essential closed curve. Moreover, since a weaving motif contains at least two non-parallel components, this means that there exist at least two distinct characteristic loops of $\mathcal{P}$ for $(\Lambda,L)$ which are isotopic to non-parallel threads in $T^2$. Similarly, if $\big(\mathcal{P},(\Lambda,L) \big)$ creates a polycatenane motif, then since all its component are null-homotopic in $T^2$, each characteristic loop of $\mathcal{P}$ for $(\Lambda,L)$ must be correspond to null-homotopic closed curve. The case of mixed motif follows directly as it contains at least one essential component and one null-homotopic component. Finally, recall that the lifts to $\mathbb{E}^2$ of weaving, polycatenane and mixed motifs do not contain components with self-crossings. This means that if a curve isotopic to a characteristic loop of $\mathcal{P}$ for $(\Lambda,L)$ contains a self-intersection point, then this point divides this curve into two non-trivial divided curves. 

Conversely, to prove the reverse implications, consider now a set of characteristic loops in $\mathcal{P}$ which does not admit trivial divided curves. Then, if this set contains only null-homotopic closed curves, it follows that $\big(\mathcal{P},(\Lambda,L) \big)$ is the scaffold of a polycatenane. 
Otherwise, if it contains at least a null-homotopic and an essential closed curve, then $\big(\mathcal{P},(\Lambda,L) \big)$ is the scaffold of a mixed motif. Finally, if the set contains only essential closed curves, such that two of them are not parallels, then $\big(\mathcal{P},(\Lambda,L) \big)$ is the scaffold of a weaving motif.
\end{proof}

\begin{figure}[ht]
\centering
   \includegraphics[width=5in]{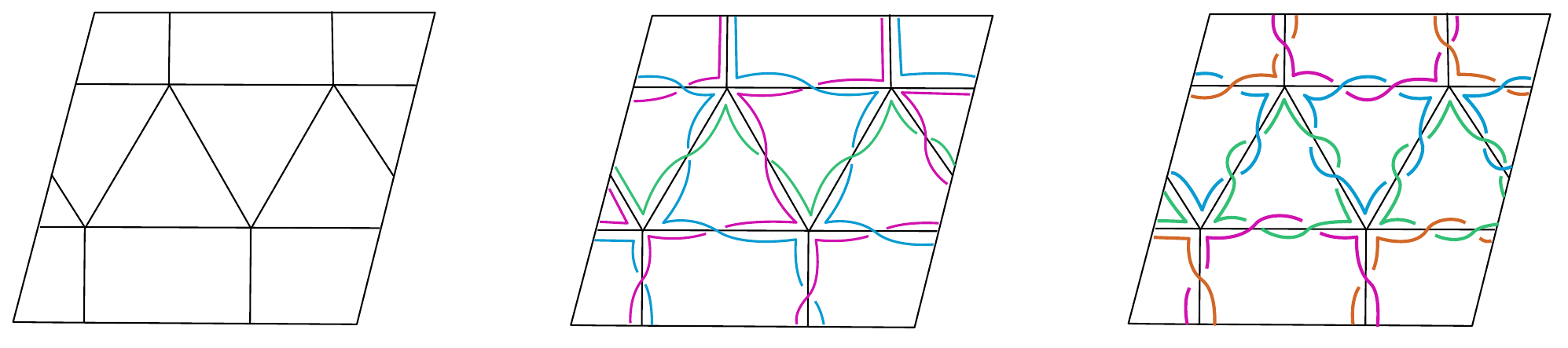}
      \caption{\label{degree5} 
      On the left, a $\mathcal{T}$-cell that is transformed by the method $(\Lambda,L) = (B_r , 1)$ to form a weaving motif (in the middle), and by the method $(\Lambda,L) = (B_r , 2)$ to form a polycatenane motif (on the right).}
\end{figure}

Examples illustrating Proposition~\ref{prop:characteristic} are shown in Figure~\ref{degree5}. Indeed by walking on the different curves generated by the polygonal link methods used as examples, we observe that the curve components of a weaving motif (in the middle) satisfy the conditions of Proposition~\ref{prop:characteristic} (1).Conversely, the curve components of a polycatenane motif (on the right) satisfy the conditions of Proposition~\ref{prop:characteristic} (2).

To investigate whether a characteristic loop in $T^2$ satisfies Proposition~\ref{prop:characteristic}, we use its polygonal description as outlined in Subsections~\ref{sec:3.2.1}, \ref{sec:3.2.2}, and \ref{sec:3.2.3}. Our approach now considers an algebraic characterization, where we use the fundamental group $\pi_1(T^2)$ to distinguish the homotopy types of characteristic loops. In particular, the fundamental group $\pi_1(T^2)$ is generated by two elements, $\alpha$ and $\beta$, corresponding to the meridian and longitude of the torus, respectively. We can also identify the flat torus with a parallelogram whose opposite sides are labeled $\alpha$ and $\beta$, accordingly (see Figure~\ref{studycase} for an illustration). 

A characteristic loop in $T^2$ can then be described as a \textit{reduced cyclic torus word} in $\pi_1(T^2)$, where a reduced word eliminates any consecutive inverse pairs of the generators, such as $g g^{-1}$. Furthermore, it is well-known that in the commutative group $\pi_1(T^2)$, a null-homotopic closed curve in $T^2$ is characterized by the trivial word $\alpha \beta \alpha^{-1} \beta^{-1}$, while an essential curve is described by a non-trivial torus word, as detailed in Chapter 1 of \cite{FarbMargalit+2012}.

\smallbreak

Let $\mathcal{P}$ be a given $\mathcal{T}$-cell and $(\Lambda,L)$ be one of the three polygonal link methods. Let also $S= (e_1^+, e_2^+, \cdots , e_l^+)$ denote the the set of oriented edges of $\mathcal{P}$, and let $F_S$ be the free group over $S$. 

Then, for each characteristic loop $\Delta^{\gamma}_{(\Lambda,L)} = (\delta_0^{\pm}, \delta_1^{\pm}, \cdots , \delta_k^{\pm})_{(\Lambda,L)}$, where for all $i \in \{1, \cdots, n\}$, $\delta_i \in \{e_1, \cdots, e_n\}$, we assign an \textit{edge word} in $F_S$: 
$$w(\Delta^{\gamma}_{(\Lambda,L)}) = \delta_0^{\pm}.\delta_1^{\pm}.\cdots .\delta_k^{\pm},$$ 
where each $\delta_i$ is an oriented edge in $\mathcal{P}$ and the $\pm$ sign reflects whether the orientation of $\delta_i$ matches the characteristic loop's direction.

In particular, we classify the edge word $w(\Delta^{\gamma}_{(\Lambda,L)})$ as follows:
\begin{itemize}  
    \item [a.] \textit{simple}: if each $\delta_i$ is distinct and if each pair of adjacent elements $(\delta_i, \delta_{i+1})$ in $w(\Delta^{\gamma}_{(\Lambda,L)})$ has a different common vertex in $\mathcal{P}$.
    \smallbreak
    \item [b.] \textit{trivial}: if the word can be written as the product of an edge word $g$ and its inverse, i.e., $w(\Delta^{\gamma}_{(\Lambda,L)}) = g.g^{-1}$, corresponding to a null-homotopic loop.
    \smallbreak
    \item [c.] \textit{knotted}, if there exists a subword $g$ such that $w(\Delta^{\gamma}_{(\Lambda,L)})=\delta_0^{\pm}.\cdots.\delta_i^{\pm}.g.\delta_j^{\pm}.\cdots.\delta_k^{\pm}$, where $\delta_i$ and $\delta_j$ share a common vertex of degree four. In this case, the subword $g$ is said to \textit{divide} the characteristic loop.
\end{itemize}

\bigbreak
\noindent We illustrate this observation with the following example.

\begin{example}
Let $\mathcal{P}$ be a $\mathcal{T}$-cell of an hexagonal lattice to which we arbitrarily label and orient its edges as illustrated in Figure~\ref{studycase}. 

\begin{figure}[ht]
\centering
   \includegraphics[width=2.5in]{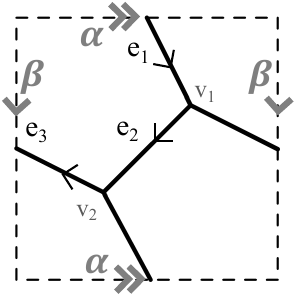}
      \caption{\label{studycase} A $\mathcal{T}$-cell of a doubly periodic hexagonal tiling.}
\end{figure}

\begin{itemize}  

  \item [1.] To predict the type of motif generated when applying $(\Lambda,L) = (\mathit{C}_r,2)$ to $\mathcal{P}$, we first enumerate the characteristic loops covered by each strand after transformation, using the polygonal chain description provided in subsection~\ref{sec:3.2.2}. In this example, all the characteristic loops are equivalent, up to cyclic permutation and reversing of all orientations, to the following polygonal chain, 
  $$\Delta^{\gamma}_{(\mathit{C}_r,2)} = (e_1^+, e_3^-, e_2^-, e_1^-, e_3^+, e_2^+)$$
  Here, the symbols $\pm$ denote the orientation of the edges, as illustrated in Figure~\ref{studycase}, with a positive sign indicating the given direction and a negative sign indicating the reversed direction. The corresponding edge word is given by, 
$$w(\Delta^{\gamma}_{(\mathit{C}_r,2)}) = e_1^+ . e_3^- . e_2^- . e_1^- . e_3^+ . e_2^+,$$
  which is trivial. Indeed, following the torus orientation, the torus word of this cycle, denoted by $w'(\Delta^{\gamma}_{(\mathit{C}_r,2)})$, can be decomposed as follows,
\begin{align*}   
w'(e_3^- . e_1^+) &= \alpha . \beta \\
w'(e_3^- . e_2^- . e_1^-) &= \alpha . \beta ^{-1} \\
w'(e_1^- . e_3^+) &= \beta ^{-1} . \alpha^{-1} \\
w'(e_3^+ . e_2^+ . e_1^+) &= \alpha ^{-1} . \beta.
\end{align*}
Thus,
$$w'(\Delta^{\gamma}_{(\mathit{C}_r,2)})= \alpha . \beta . \alpha . \beta ^{-1} . \beta ^{-1} . \alpha^{-1} . \alpha ^{-1} . \beta = 0 $$

By applying the same process to the other curve components constructed by the polygonal link method, we conclude that a trivial torus word can be assigned to each characteristic loop in $\mathcal{P}$ for this example. Consequently, we will demonstrate through Theorem~\ref{thm:polycatenaneprediction} that the pair $\big(\mathcal{P}, (\mathit{C}_r, 2)\big)$, where $\mathcal{P}$ is a $\mathcal{T}$-cell of a hexagonal tiling, always generates a polycatenane. It is also worth noting that the same conclusion applies to any choice of an even number of twists $m$ in the polygonal link method.

\smallbreak

  \item [2.] Next, to predict the type of motif created when applying $(\Lambda,L) = (\mathit{B}_r,3)$ to $\mathcal{P}$, we begin by listing the characteristic loops covered by each strand after transformation, using the polygonal chain description provided in subsection~\ref{sec:3.2.3}. Here, all characteristic loops can be described by the following simple edge words, up to cyclic permutation and reversing all orientations,
\begin{align*}   
w(\Delta^1_{(\mathit{B}_r,3)}) &= e_1^+ . e_3^- , \\
w(\Delta^2_{(\mathit{B}_r,3)}) &= e_1^+ . e_2^+ , \\
w(\Delta^3_{(\mathit{B}_r,3)}) &= e_3^+ . e_2^+ .
\end{align*}

Moreover, these characteristic loops are closed curves with non-trivial and distinct torus words, namely $\alpha . \beta$, $\beta$, and $\alpha^{-1}$, respectively. Thus, we will prove by Theorem~\ref{thm:weavingprediction} that the pair $\big(\mathcal{P},(\mathit{B}_r,3)\big)$, where $\mathcal{P}$ is a $\mathcal{T}$-cell of a hexagonal tiling, always generates a twisted weave. Note that the same conclusion also holds for any odd number of twists $m$.
 \end{itemize} 
\end{example}

The above classification of edge words allows us to algebraically characterize the types of motifs generated by the polygonal link methods. We are now in a position to predict if a specific type of motif can be generated from a given $\mathcal{T}$-cell and a chosen polygonal link method $(\Lambda, L)$, as stated in the following main theorems.

\begin{theorem}\label{thm:weavingprediction}\textbf{(Weaving Motif Prediction)} 
Let $\mathcal{P}$ be a $\mathcal{T}$-cell and $(\Lambda,L)$ a polygonal link method. Then, $\big(\mathcal{P},(\Lambda,L) \big)$ generates a weaving motif if and only if every characteristic loop $\Delta^{\gamma}_{(\Lambda,L)}$ in $\mathcal{P}$ can be written as a simple or knotted edge word $w$ satisfying all the following conditions:
  \begin{itemize}
      \item [i.] the corresponding torus word of $w$ is non-trivial in $\pi_1(T^2)$, 
      \smallbreak
      \item [ii.] there exist at least two distinct characteristic loops whose corresponding reduced torus words are not homotopically equivalent in $\pi_1(T^2)$,
      \smallbreak
      \item [iii.] if $w$ is knotted, then any dividing subword must have a non-trivial corresponding torus word in $\pi_1(T^2)$.
  \end{itemize}
\end{theorem}

\begin{proof} 
The forward implication follows from the definitions of untwisted and twisted weaving motifs (Definitions~\ref{def:untwistedweave}, \ref{def:twistedweave}). 
Let $\big(\mathcal{P}, (\Lambda,L)\big)$ generate a weaving motif. A weaving motif consists of essential closed curves in $T^2$, so their corresponding torus words in $\pi_1(T^2)$ must be non-trivial. Therefore, by definition, every characteristic loop in a weaving motif corresponds to a non-trivial torus word, satisfying condition (i). 

Moreover, by definition, a weaving motif must contain at least two distinct, non-parallel components. In terms of homotopy terms, this implies that there must be at least two characteristic loops that represent distinct homotopy classes in $T^2$. Algebraically, this means that their corresponding torus words in $\pi_1(T^2)$ are not equivalent, satisfying condition (ii).

Finally, consider the case where a characteristic loop is knotted, meaning it has a self-crossing. By (1) of Proposition~\ref{prop:characteristic}, this loop can be divided into two divided curves and for the characteristic loop to remain essential, both divided curves must be essential, meaning that they both have non-trivial torus words. Therefore, if $w$ is knotted, any dividing subword must have a non-trivial corresponding torus word in $\pi_1(T^2)$, satisfying condition (iii).

To prove the reverse implication, we assume that the set of characteristic loops of $\mathcal{P}$ for $(\Lambda,L)$ satisfies conditions (i), (ii), and (iii). We will show that $\big(\mathcal{P}, (\Lambda,L)\big)$ generates a weaving motif. By condition (i), each characteristic loop corresponds to a non-trivial torus word, meaning that all characteristic loops are essential closed curves. 
Next, by condition (ii), there are at least two distinct characteristic loops whose corresponding torus words are not equivalent in $\pi_1(T^2)$. Two loops with non-equivalent torus words represent distinct homotopy classes, thus the loops are \textit{not parallel} in $T^2$.
Finally, consider the case where a characteristic loop is knotted. By condition (iii), if an essential characteristic loop is knotted, its dividing subwords must correspond to essential divided curves, which also have non-trivial torus words.
Therefore, $\big(\mathcal{P}, (\Lambda, L) \big)$ generates a weaving motif, since all characteristic loops are essential, and the motif contains at least two non-parallel components. This completes the proof of the theorem.
\end{proof}

\begin{theorem}\label{thm:polycatenaneprediction}\textbf{(Polycatenane Motif Prediction)} 
Let $\mathcal{P}$ be a $\mathcal{T}$-cell and $(\Lambda,L)$ a polygonal link method. Then, $\big(\mathcal{P},(\Lambda,L) \big)$ generates a polycatenane motif if and only if every characteristic loop $\Delta^{\gamma}_{(\Lambda,L)}$ in $\mathcal{P}$ can be written as an edge word $w$ satisfying all the following conditions:
  \begin{itemize}
      \item [i.] the corresponding torus word of $w$ is trivial in $\pi_1(T^2)$,
      \smallbreak
      \item [ii.] if $w$ is knotted, then any of its dividing subwords has a non-trivial corresponding torus word in $\pi_1(T^2)$.
  \end{itemize}
\end{theorem}

\begin{proof}
We begin with the forward implication. 

Let $\big(\mathcal{P}, (\Lambda,L)\big)$ generate a polycatenane motif. By Definition~\ref{def:polycatenane}, a polycatenane consists of null-homotopic closed curves in $T^2$. Since null-homotopic curves in $T^2$ correspond to the trivial element in $\pi_1(T^2)$, this implies that the torus word associated with each characteristic loop is trivial, satisfying condition (i). Condition (ii) is satisfied using the same argument as in the proof of Theorem~\ref{thm:weavingprediction}.

To prove the reverse implication, assume that the set of characteristic loops of $\mathcal{P}$ for $(\Lambda,L)$ satisfies conditions (i) and (ii). We will show that $\big(\mathcal{P}, (\Lambda,L)\big)$ generates a polycatenane motif.
By condition (i), the torus word associated with each characteristic loop is trivial, which implies that all characteristic loops are null-homotopic closed curves in $T^2$, as needed for a polycatenane motif by (2) of Proposition~\ref{prop:characteristic}. Then, once again, condition (ii) is treated as in the proof of Theorem~\ref{thm:weavingprediction}.
Therefore, all characteristic loops are null-homotopic closed curves in $T^2$, and the motif contains only null-homotopic components. By the definition of a polycatenane motif, $\big(\mathcal{P}, (\Lambda,L)\big)$ generates a polycatenane motif. This completes the proof of the theorem.
\end{proof}

Finally, the corollary stated below follows directly from (3) of Proposition~\ref{prop:characteristic}, Theorems~\ref{thm:weavingprediction} and~\ref{thm:polycatenaneprediction}, as a mixed motif contains both essential and null-homotopic components by Definition~\ref{def:mixed}. Thus, a detailed proof is omitted to avoid redundancy.

\begin{corollary}\label{thm:mixedprediction}\textbf{(Mixed Motif Prediction)} 
Let $\mathcal{P}$ be a $\mathcal{T}$-cell and $(\Lambda,L)$ be a polygonal link method. Then, $\big(\mathcal{P},(\Lambda,L) \big)$ generates a mixed motif if and only if the set of all the characteristic loops $\Delta_{(\Lambda,L)}$ in $\mathcal{P}$ contains at least a component with corresponding edge word satisfying the conditions (i) and (iii) of Theorem~\ref{thm:weavingprediction}, and a component with corresponding edge word satisfying Theorem~\ref{thm:polycatenaneprediction}.
\end{corollary}

%


\end{document}